\documentclass[12pt]{article}
\usepackage{amsmath,amsthm,amsfonts,amssymb,amscd}
\def\qed{{\unskip\nobreak\hfil\penalty50
\hskip2em\hbox{}\nobreak\hfil$\square$
\parfillskip=0pt \finalhyphendemerits=0\par}\medskip}

\def\Aut{{\mathrm {Aut}}}
\def\End{{\mathrm {End}}}

\def\Tr{{\mathrm {Tr}}}

\def\a{\alpha}
\def\b{\beta}

\def\l{\lambda}

\def\phi{\varphi}

\def\th{\theta}

\def\Om{\Omega}

\def\p{{\pi}}

\def\emptyset{\varnothing}
\def\setminus{\smallsetminus}

\def\Diff{{\mathrm {Diff}}}

\def\PSL{PSU(1,1)}
\def\supp{{\rm supp}}

\def\Aut{{\mathrm {Aut}}}
\def\End{{\mathrm {End}}}

\def\Tr{{\mathrm {Tr}}}

\def\a{\alpha}
\def\b{\beta}

\def\l{\lambda}
\def\phi{\varphi}

\def\th{\theta}

\def\Om{\Omega}

\def\p{{\pi}}

\newtheorem{theorem}{Theorem}[section]
\newtheorem{lemma}[theorem]{Lemma}
\newtheorem{conjecture}[theorem]{Conjecture}
\newtheorem{corollary}[theorem]{Corollary}
\newtheorem{definition}[theorem]{Definition}

\newtheorem{proposition}[theorem]{Proposition}
\newtheorem{remark}[theorem]{Remark}
\newtheorem{example}[theorem]{Example}

\def\emptyset{\varnothing}
\def\setminus{\smallsetminus}

\def\Diff{{\mathrm {Diff}}}

\def\PSL{PSU(1,1)}

\def\res{\!\restriction\!}
\def\A{{\cal A}}

\def\B{{\cal B}}

\def\I{{\cal I}}

\def\L{{\cal L}}

\def\H{{\cal H}}

\def\Z{{\mathbb Z}}
\def\R{{\mathbb R}}
\def\Co{{\mathbb C}}
\def\1#1{{\bf #1}}
\def\2#1{{\mathcal #1}}
\def\3#1{{\sl #1}}
\def\4#1{{\tt #1}}
\def\5#1{{\sf #1}}
\def\6#1{{\mathfrak #1}}
\def\7#1{{\mathbb #1}}
\def\<{{\langle}}
\def\>{{\rangle}}

\newcommand{\pf}{{\it Proof. }}
\renewcommand{\qed}{\ \hfill $\blacksquare$}

\newcommand{\bdefin}{\begin{definition}}
\newcommand{\blemma}{\begin{lemma}}
\newcommand{\bprop}{\begin{proposition}}
\newcommand{\btheor}{\begin{theorem}}
\newcommand{\bcoro}{\begin{corollary}}
\newcommand{\bconj}{\begin{conjecture}}
\newcommand{\edefin}{\end{definition}}
\newcommand{\elemma}{\end{lemma}}
\newcommand{\eprop}{\end{proposition}}
\newcommand{\etheor}{\end{theorem}}
\newcommand{\ecoro}{\end{corollary}}
\newcommand{\econj}{\end{conjecture}}
\newcommand{\brem}{\begin{remark}}
\newcommand{\erem}{\end{remark}}

\newcommand{\ba}{\begin{array}}
\newcommand{\ea}{\end{array}}
\newcommand{\bea}{\begin{eqnarray}}
\newcommand{\eea}{\end{eqnarray}}
\newcommand{\bean}{\begin{eqnarray*}}
\newcommand{\eean}{\end{eqnarray*}}

\newcommand{\inz}{\epsilon}

\title{\huge Conformal Nets Associated With Lattices And Their Orbifolds\\}
\author{
{\sc Chongying Dong}\footnote{Supported in part by NSF and
a research grant from Univercity of California at Santa Cruz.}\\
Department of Mathematics\\
University of California at Santa Cruz\\
Santa Cruz, CA 95064\\
E-mail: {\tt dong@math.ucsc.edu}\\
{}\\
{\sc Feng Xu}\footnote{Supported in part by NSF.}\\
Department of Mathematics\\
University of California at Riverside\\
Riverside, CA 92521\\
E-mail: {\tt xufeng@math.ucr.edu}}
\begin{document}
\date{}
\maketitle

\begin{abstract}
In this paper we study representations of conformal nets associated with 
positive definite even lattices and their orbifolds with respect to isometries
of the lattices. Using previous general results on orbifolds, we give a 
list of all irreducible representations of the orbifolds, which generate
a unitary modular tensor category. 2000MSC:81R15, 17B69.
\end{abstract}
\newpage

\section{Introduction}

Let $\A$ be a completely rational conformal net (cf. \S2 and (\ref{abr})
following \cite{KLM}).   
Let $\Gamma$ be a finite group acting properly on $\A$ (cf. definition
(\ref{p'})). The starting point of this paper is Th. \ref{orb} proved
in \cite{Xu3} which states that the fixed point subnet (the orbifold)
$\A^\Gamma$ is also  completely rational, and by \cite{KLM} 
$\A^\Gamma$ has finitely many irreducible representations which are 
divided into two classes: the ones that are obtained from the 
restrictions of a representation of $\A$ to $\A^\Gamma$ which 
are called untwisted representations, and the ones which are
twisted. It follows from 
 Th. \ref{orb} that twisted representation of $\A^\Gamma$
always exists if  $\A^\Gamma\neq \A.$ The motivating question for
this paper is how to construct these   twisted representations of 
$\A^\Gamma.$\par
It turns out that all  representations of 
$\A^\Gamma $ are closely related to the solitons of $\A.$ 
In \cite{KLX}, we construct solitons for 
the affine and permutation orbifolds. In this paper, we give a construction
of solitons for the case of conformal nets associated with positive definite
even lattices and isomteries of the lattices. We note that 
conformal nets associated 
with special lattices have appeared before in \S3 of \cite{Xu1} and more
recently in \cite{KL}. 
Our solitons correspond to
``twisted representations'' of the corresponding vertex operator algebras
(VOAs), and such twisted representations have been studied in
\cite{L}, \cite{FLM}, \cite{FLM1}, \cite{D1}, \cite{DL} and more recently \cite{BK}. We will show that these
twisted representations in the VOA sense indeed give rise to solitons 
(cf. Definition \ref{generalsoliton} and Prop. \ref{generalnormal}). Compared
to the constructions of \cite{KLX}, the notable difference is that the
twisted representations are not related to the untwisted representations
in a simple way as in the  affine and permutation orbifolds: in the later
cases twisted representations are constructed on the same spaces as that of
untwisted representations, only with a ``twisted'' action, while in the
cases of lattices the spaces are different. Hence it is nontrivial to show
that   in the cases of lattices these
twisted representations in the VOA sense indeed give rise to solitons. Such 
questions were already encountered in \cite{Xu2} when the lattice has
rank one, and were solved by identifying the orbifolds as special cosets
in \cite{Xu1}. However this method of  \cite{Xu2} does not work for
general lattices. Our solution to this question consists of two steps.
First we show that if the inner product on the lattice has some integral 
property, the question can be solved (cf. Prop. \ref{Normality1}) 
by using a covering homorphism (cf. Prop. \ref{N-fold}). For general lattice
$Q,$ we choose a sublattice $P\subset Q$ with finite index 
such that the restriction of the inner product from $Q$ to $P$ has the desired
integral property as in the first step. By exploiting the related 
group structures (cf. Lemma \ref{homo}) and results from the first step, 
we then give Definition \ref{generalsoliton} and show that they have 
the right properties in Prop. \ref{generalnormal}. \par
To show that these solitons give rise to all irreducible representations
of the orbifolds, we use the same strategy as in \cite{KLX} which is to
compute the index of solitons and use Th. \ref{orb}.  Here we make use of
the large group of (local) automorphisms (cf. Definition \ref{auto A_Q}) of 
the conformal nets associated with the
lattices, and an exhaustion trick 
(cf. the paragraph before Th. \ref{allirrep}) to prove   
Th. \ref{allirrep}. Th. \ref{allirrep} is the main result of this paper:
it gives a list of all irreducible representations of the orbifold from
solitons, and they generate a unitary modular tensor category 
(cf. \cite{Tu}). We expect that this result will have
applications in many concrete examples, and we plan to address such 
applications in the future.\par
The rest of the paper is organized briefly as follows: after introducing
basic notions such as conformal nets, their representations, complete
rationality, orbifolds and solitons in \S2, we consider conformal nets
associated with positive definite even lattices in \S3. The main result
in \S3 is Cor. \ref{A_Q rational}. In \S4 we consider the constructions of
solitons. As mentioned above we do this in two steps: the first construction
in a special case is given in Definition \ref{basic soliton} and its
properties are studied in Prop. \ref{Normality1}, Prop. \ref{Irred}, 
Prop. \ref{Irred2}, Prop. \ref{sigma_1}; the general case is considered in
subsection \ref{gs} where   Th. \ref{allirrep} is proved.

\section{Preliminaries}
In this section we review some basic concepts of conformal nets
which will be used. See \S2 and \S3 of \cite{KLX} for a more detailed
review. 
\subsection{Conformal nets}

We denote by $\I$ the family of proper open intervals of $S^1$. 
A {\it net} $\A$ of von Neumann algebras on $S^1$ is a map 
\[
I\in\I\to\A(I)\subset B(\H)
\]
from $\I$ to von Neumann algebras on a fixed Hilbert space $\H$
that satisfies:
\begin{itemize}
\item[{\bf A.}] {\it Isotony}. If $I_{1}\subset I_{2}$ belong to 
$\I$, then
\begin{equation*}
 \A(I_{1})\subset\A(I_{2}).
\end{equation*}
\end{itemize}
If $E\subset S^1$ is any region, we shall put 
$\A(E)\equiv\bigvee_{E\supset I\in\I}\A(I)$ with $\A(E)=\mathbb C$ 
if $E$ has empty interior (the symbol $\vee$ denotes the von Neumann 
algebra generated). 

The net $\A$ is called {\it local} if it satisfies:
\begin{itemize}
\item[{\bf B.}] {\it Locality}. If $I_{1},I_{2}\in\I$ and $I_1\cap 
I_2=\emptyset$ then 
\begin{equation*}
 [\A(I_{1}),\A(I_{2})]=\{0\},
 \end{equation*}
where brackets denote the commutator.
\end{itemize}
The net $\A$ is called {\it M\"{o}bius covariant} if in addition 
satisfies
the following properties {\bf C,D,E,F}:
\begin{itemize}
\item[{\bf C.}] {\it M\"{o}bius covariance}. 
There exists a strongly 
continuous unitary representation $U$ of the M\"{o}bius group 
$\bold G$ on $\H$ such that
\begin{equation*}
 U(g)\A(I) U(g)^*\ =\ \A(gI),\quad g\in {\bold G},\ I\in\I.
\end{equation*}
Note that this implies $\A(\bar I)=\A(I)$, $I\in\I$ (consider a 
sequence 
of elements $g_n\in{\bold G}$ converging to the identity such that 
$g_n\bar I \nearrow I$).
\item[{\bf D.}] {\it Positivity of the energy}. The generator of the 
one-parameter rotation subgroup of $U$ (conformal Hamiltonian) is 
positive. 
\item[{\bf E.}] {\it Existence of the vacuum}. There exists a unit 
$U$-invariant vector $\Omega\in\H$ (vacuum vector), and $\Omega$ 
is cyclic for the von Neumann algebra $\bigvee_{I\in\I}\A(I)$.
\end{itemize}
By the Reeh-Schlieder theorem $\Omega$ is cyclic and separating for 
every fixed $\A(I)$. The modular objects associated with 
$(\A(I),\Omega)$ have a geometric meaning
\[
\Delta^{it}_I = U(\Lambda_I(2\pi t),\qquad J_I = U(r_I)\ .
\]
Here $\Lambda_I$ is a canonical one-parameter subgroup of ${\bold G}$ and $U(r_I)$ is a 
antiunitary acting geometrically on $\A$ as a reflection $r_I$ on $S^1$. 

This imply {\em Haag duality}: 
\[
\A(I)'=\A(I'),\quad I\in\I\ ,
\]
where $I'$ is the interior of $S^1\setminus I$.

\begin{itemize}
\item[{\bf F.}] {\it Irreducibility}. $\bigvee_{I\in\I}\A(I)=B(\H)$. 
Indeed $\A$ is irreducible iff
$\Om$ is the unique $U$-invariant vector (up to scalar multiples). 
Also  $\A$ is irreducible
iff the local von Neumann 
algebras $\A(I)$ are factors. In this case they are III$_1$-factors in 
Connes classification of type III factors
(unless $\A(I)=\mathbb C$ identically).
\end{itemize}
By a {\it conformal net} (or diffeomorphism covariant net)  
$\A$ we shall mean a M\"{o}bius covariant net such that the following 
holds:
\begin{itemize}
\item[{\bf G.}] {\it Conformal covariance}. There exists a projective 
unitary representation $U$ of $\Diff(S^1)$ on $\H$ extending the unitary 
representation of $\PSL$ such that for all $I\in\I$ we have
\begin{gather*}
 U(g)\A(I) U(g)^*\ =\ \A(gI),\quad  g\in\Diff(S^1), \\
 U(g)xU(g)^*\ =\ x,\quad x\in\A(I),\ g\in\Diff(I'),
\end{gather*}
\end{itemize}
where $\Diff(S^1)$ denotes the group of smooth, positively oriented 
diffeomorphism of $S^1$ and $\Diff(I)$ the subgroup of 
diffeomorphisms $g$ such that $g(z)=z$ for all $z\in I'$.
A (DHR) representation $\pi$ of $\A$ on a Hilbert space $\H$ is a map 
$I\in\I\mapsto  \pi_I$ that associates to each $I$ a normal 
representation of $\A(I)$ on $B(\H)$ such that
\[
\pi_{\tilde I}\res\A(I)=\pi_I,\quad I\subset\tilde I, \quad 
I,\tilde I\subset\I\ .
\]
$\pi$ is said to be M\"obius (resp. diffeomorphism) covariant with
positive energy if 
there is a projective unitary representation $U_{\pi}$ of ${\bold G}$ (resp. 
$\Diff^{(\infty)}(S^1)$, the infinite cover of $\Diff(S^1)$ ) 
with positive energy (the generator of the rotation subgroup
$S^1$ has non-negative spectrum) on $\H$ such that
\[
\pi_{gI}(U(g)xU(g)^*) =U_{\pi}(g)\pi_{I}(x)U_{\pi}(g)^*
\]
for all $I\in\I$, $x\in\A(I)$ and $g\in {\bold G}$ (resp. 
$g\in\Diff^{(\infty)}(S^1)$). Note that if $\pi$ is irreducible and 
diffeomorphism covariant then $U$ is indeed a projective unitary 
representation of $\Diff(S^1)$.\par
\subsection{The orbifolds}
Let ${\A}$ be an irreducible conformal net on a Hilbert space
${\H}$ and let $\Gamma$ be a finite group. Let $V:\Gamma\rightarrow U({\H})$
be a  unitary representation of $\Gamma$ on ${\H}$. 
If $V:\Gamma\rightarrow U({\H})$ is not faithful, we set $\Gamma':= \Gamma
/{\rm ker} V$.
\begin{definition} \label{p'}
We say that $\Gamma$ acts properly on ${\A}$ if the following conditions
are satisfied:\par
(1) For each fixed interval $I$ and each $g\in \Gamma$, 
$\alpha_g (a):=V(g)aV(g^*) \in {\A}(I), \forall a\in
{\A}(I)$; \par
(2) For each  $g\in \Gamma$, $V(g)\Omega = \Omega, \forall g\in \Gamma$.\par
\label{Definition 2.1}
\end{definition}
We note that if $\Gamma$ acts properly, then $\Gamma$ commutes with
${\mathbf G}$. \par
Define ${\A}^\Gamma(I):={\A}(I)P_0$ on ${\H}_0$ where $\H_0:=\{ x\in \H| 
V(g)x=x, \forall g\in \Gamma \}$ and $P_0$ is the projection
from $\H$ to $\H_0.$  The unitary
representation $U$ of ${\mathbf G}$ on ${\H}$ restricts to
an  unitary
representation (still denoted by $U$) of ${\mathbf G}$ on ${\H}_0$.
Then:
\bprop 
The map $I\in {\I}\rightarrow {\A}^G(I)$ on $ {\H}_0$ 
together with the  unitary
representation (still denoted by $U$) of ${\mathbf G}$ on ${\H}_0$
is an
irreducible M\"{o}bius covariant net.
\label{Prop.2.1}
\eprop
The irreducible  M\"{o}bius covariant net in Prop.  \ref{Prop.2.1} 
will be denoted by
${\A}^\Gamma$ and will be called the {\it orbifold of ${\A}$}
with respect to $\Gamma$. We note that by definition 
${\A}^\Gamma= {\A}^{\Gamma'}$. \par

\subsection{Complete rationality}
By an interval of the circle we mean an open connected
proper subset of the circle. If $I$ is such an interval then
$I'$ will denote the interior of the complement of $I$ in the circle.
We will denote by ${\I}$ the set of such intervals.
Let $I_1, I_2\in {\I}$. We say that $I_1, I_2$ are disjoint if
$\bar I_1\cap \bar I_2=\emptyset$, where $\bar I$
is the closure of $I$ in $S^1.$
Denote by ${\I}_2$ the set of unions of disjoint 2 elements
in ${\I}$. Let ${\A}$ be an irreducible conformal net
as in \S2.1. For $E=I_1\cup I_2\in{\I}_2$, let
$I_3\cup I_4$ be the interior of the complement of $I_1\cup I_2$ in 
$S^1$ where $I_3, I_4$ are disjoint intervals. 
Let 
$$
{\A}(E):= A(I_1)\vee A(I_2), 
\hat {\A}(E):= (A(I_3)\vee A(I_4))'.
$$ Note that ${\A}(E) \subset \hat {\A}(E)$.
Recall that a net ${\A}$ is {\it split} if ${\A}(I_1)\vee
{\A}(I_2)$ is naturally isomorphic to the tensor product of
von Neumann algebras ${\A}(I_1)\otimes
{\A}(I_2)$ for any disjoint intervals $I_1, I_2\in {\I}$.
${\A}$ is {\it strongly additive} if ${\A}(I_1)\vee
{\A}(I_2)= {\A}(I)$ where $I_1\cup I_2$ is obtained
by removing an interior point from $I$.
\bdefin\label{abr}
\cite{KLM}
${\A}$ is said to be completely  rational, or $\mu$-rational, if
${\A}$ is split, strongly additive, and 
the index $[\hat {\A}(E): {\A}(E)]$ is finite for some
$E\in {\I}_2$ . The value of the index
$[\hat {\A}(E): {\A}(E)]$ (it is independent of 
$E$ by Prop. 5 of \cite{KLM}) is denoted by $\mu_{{\A}}$
and is called the $\mu$-index of ${\A}$. If 
the index $[\hat {\A}(E): {\A}(E)]$ is infinity for some
$E\in {\I}_2$, we define the $\mu$-index of ${\A}$ to be
infinity.
\label{Definition 2.2}
\edefin
The following theorem is proved in \cite{Xu3}:
\btheor\label{orb}
Let ${\A}$ be an irreducible conformal net and let $\Gamma$
be a finite group acting properly on ${\A}$. Suppose that 
${\A}$ is completely rational or $\mu$-rational as in definition
2.2. Then:\par
(1): ${\A}^\Gamma$ is completely rational or $\mu$-rational and  
$\mu_{{\A}^\Gamma}= |\Gamma'|^2 \mu_{{\A}}$; \par
(2): There are only a finite number of irreducible covariant 
representations of ${\A}^\Gamma$, and they give rise to a unitary modular
category as defined in II.5 of \cite{Tu} by the construction as given in
\S1.7 of \cite{Xu4}.
\label{Th.2.6}
\etheor
Suppose that ${\A}$ and $\Gamma$ satisfy the assumptions of Th.\ref{orb}. 
Then ${\A}^\Gamma$ has only a  
finite number of irreducible representations $\dot\lambda$ and
$$
\sum_{\dot\lambda}d(\lambda)^2 = \mu_{{\A}^\Gamma}= |\Gamma'|^2 \mu_{{\A}} 
$$
where we use $d(\lambda)$ to denote the statistical dimension or
the square root of index (cf. \cite{J} and \cite{PP}).
\subsection{Solitons}
\label{solitonsection}
Let $\xi\in S^1,$ and identify $\R$ with $S^1\setminus \{\xi\}\simeq (0,1).$
Denote by $\I_0$ the set of open, 
connected, non-empty, proper subsets of $\mathbb R$, thus $I\in\I_0$ 
iff $I$ is an open interval or half-line (by an interval of $\mathbb 
R$ we shall always mean a non-empty open bounded interval of $\mathbb 
R$).

Given a net $\A$ on $S^1$ we shall denote by $\A_0$ its restriction 
to $\mathbb R = S^1\setminus\{-1\}$. Thus $\A_0$ is an isotone map on
$\I_0$, that we call a \emph{net on $\mathbb R$}.

A representation $\pi$ of $\A_0$ on a Hilbert space $\H$ is a map 
$I\in\I_0\mapsto\pi_I$ that associates to each $I\in\I_0$ a normal 
representation of $\A(I)$ on $B(\H)$ such that
\[
\pi_{\tilde I}\res\A(I)=\pi_I,\quad I\subset\tilde I, \quad 
I,\tilde I\in\I_0\ .
\]
A representation $\pi$ of $\A_0$ is also called a
\emph{soliton}. If we wish to emphasize on the dependence of $\xi\in S^1,$
we will write $\pi$ as $\pi^{(\xi)}.$
\section{Conformal nets associated with a lattice and their representations}

Let $Q$ be a positive definite even lattice. That is, $Q$ is a free
abelian group of finite rank with a positive definite $\Z$-valued 
bilinear form $\<\cdot\>$ such that $\<\a,\a\>\in 2\Z$ for all
$\a\in Q.$  $Q^*=\{\b\in \R Q|\<\b,Q\>\subset \Z\}$ is the dual lattice
of $Q.$ There exists a bimultiplicative function 
$\epsilon: Q\times Q\rightarrow \{\pm 1\}$ satisfying 
$\epsilon(\alpha,\alpha)= (-1)^{\langle \alpha,\alpha\rangle/2}, \alpha\in Q.$
Then by  bimultiplicativity 
$\epsilon(\alpha,\beta) \epsilon(\beta,\alpha)= 
(-1)^{\langle \alpha,\beta\rangle}  
, \alpha, \beta\in Q.$ Note that such 2-cocycle $\epsilon$ is unique up
to equivalence. We say  $\epsilon$ is trivial if 
$\epsilon(\alpha,\beta)=1, \forall  \alpha,\beta \in Q.$
We note that one can always choose a trivial  2-cocycle (in the 
equivalence class) if $\langle \alpha,\beta\rangle\in 4\Z, \forall  \alpha,\beta \in Q.$ 
Let $T=\R Q/Q$ be the  torus.  We will represent elements of $T$ by $e^{2\pi i h}$ for $h\in \R Q.$ Note that
$e^{2\pi i h}=1$ iff $h\in Q.$

Denote by $LT=C^{\infty}(S^1,T).$ Every element of $LT$ can be written as $e^{2\pi i f},$ where $f=f(\th): S^1\to \R Q$ with $0\leq \th\leq 1$ and
$f(\th)=\Delta_f\th+f_0+f_1(\theta).$ Here $\Delta_f\in Q$ 
is called the ``winding number''
of $f,$ and $f_1(\th)=\sum_{n\ne 0}a_ne^{2\pi in\theta}$ for some $a_n.$ $f_0$ is 
called the``zero mode'' of $f.$ The rotation group $S^1$ acts naturally on
$LT$, and we denote the action by $R_\th.$ Define
$\int_{S^1} \langle f,g\rangle d\th:=\int_{0}^{1}
 \langle f,g\rangle d\th.$

Let $\L T=LT\times S^1,$ and define a multiplication on 
$\L T$ as follows:
\begin{definition}\label{MultiplicationofLT}
$$(e^{2\pi i f},x_1)(e^{2\pi i g},x_2)=(e^{2\pi i(f+g)},x_1x_2\epsilon(\Delta f,\Delta g)
e^{\pi i[\int\<f'|g\>d\th-\<f(1)|\Delta_g\>+\frac{1}{2}\<\Delta_f|\Delta_g\>]})$$
\end{definition}

\begin{lemma}\label{LT} $\L T$ with the above multiplication is a central 
extension of $LT.$ Moreover, the action of rotation $R_\th$ on $LT$ lifts
naturally to $\L T$ and we have 
$R_{\theta}(e^{2\pi if})R_{\theta}(e^{2\pi ig})=R_{\theta}(e^{2\pi if}e^{2\pi ig})$ in $\L T$.
\end{lemma}
\pf We note that the associativity of the multiplication 
follows from the properties of 2-cocycle $\inz,$ and the rest follows by
using definitions.
\qed 
\begin{remark}
Our choice of multiplication rules in the above definition is different from that
of \cite{PS} (cf. Chapter 4 of  \cite{PS}) in the special case when $Q$ is a root 
lattice, and such a choice makes the action of rotations simpler than that
of \cite{PS}.
\end{remark} 

\begin{proposition}\label{locality} Let $f, g$ be such that $\supp e^{2\pi i f}\cap \supp e^{2\pi i g}=\emptyset$  where $\supp e^{2\pi i f}$ is defined 
to be the support of $e^{2\pi i f}$ as an element in $LT$.  Then 
$$ e^{2\pi i f} e^{2\pi i g}=
e^{2\pi i g} e^{2\pi i f}$$
as elements in $\L T.$ 
\end{proposition}

\pf Assume that $f(\theta)=\a\th+f_0+f_1(\th),$ $g(\theta)=\b\th+g_0+g_1(\th)$
and $g(0)=0$ and $g(1)=\beta.$ Then by definition \ref{MultiplicationofLT}:
$$e^{2\pi i f}e^{2\pi i g}(e^{2\pi i f})^{-1}=e^{2\pi i g}
e^{2\pi i(\frac{1}{2}\<\a,\b\>+\<\a,g_0\>-\<\b,f_0\>-\int_{S^1}\<g_1'|f_1\>d\th)}.$$
Note that
\begin{eqnarray*}
& & \int_{S^1}\<g_1'|f_1\>d\th=\int_{S^1}\<g_1'|f-\a\th-f_0\>d\th\\
& &\ \ \ =\int_{S^1}\<g'|f\>d\th-\<\a,\beta\>+\int_{S^1}\<g|\a\>d\th-\<\b|f_0\>\\
&&\ \ \ =\int_{S^1}\<g'|f\>d\th-\<\a,\beta\>+\int_{S^1}\<\b\th+g_0|\a\>d\th-\<\b|f_0\>\\
& &\ \ \ =\int_{S^1}\<g'|f\>d\th-\frac{1}{2}\<\a,\beta\>+\int_{S^1}\<g_0|\a\>d\th-\<\b|f_0\>.
\end{eqnarray*}
Since  $\supp e^{2\pi i f}\cap \supp e^{2\pi i g}=\emptyset,\int_{S^1}\<g'|f\>d\th$ is either 0 or $\<\a|\b\>.$ It follows that 
$$e^{2\pi i f}e^{2\pi i g}(e^{2\pi i f})^{-1}=e^{2\pi i g}.$$
\qed
\par
The structure of $\L T$ is well known (cf. {Page 191 of
\cite{PS}). The identity component $(\L T)^o$ of $\L T$ is canonically
a product $T\times \tilde{V}_Q,$ where $\tilde{V}_Q$ is the Heisenberg
group defined as follows: Let $W_0$ be the set of maps $f: S^1\to \R Q$ 
with winding number and zero mode being zeros. Then $\tilde{V}_Q$ is equal
to $W_0\times S^1$ as sets and multiplication is determined by
$$(f_1,\l_1)(f_2,\l_2)=(f_1+f_2,\l_1\l_2e^{\frac{i}{2}\int \<f_1'|f_2\>d\th}).$$
Let $W$ be the set of maps $f: S^1\to \R Q$ with winding number zero,
and $\tilde{W}:=\tilde{V}_Q\times \R Q.$ The following is essentially Prop. 9.5.10 in \cite{PS}:
\begin{lemma}\label{irrep of LT} We have:\par
(1)  $\tilde{V_Q}$ has a unique irreducible representation with positive
energy on a Hilbert space denoted by $S(V);$

(2) All irreducible representations of  $\tilde{W}$ with positive energy 
are of the form $S(V)_{\a}$ for $\a\in \R Q$ where $S(V)_{\a}$ is
the same as $S(V)$ as a representation of  $\tilde{V}$ and the center 
$(0,h)$ of  $\tilde{W}$ acts  on $S(V)_{\a}$ as a scalar $e^{2\pi i\<h|\a\>};$

(3) All irreducible representations of $\L T$ with positive energy are
of the form $H_{\l}=\bigoplus_{\a\in \l+Q}S(V)_\a$ for $\l\in Q^*$
where $Q^*=\{\b\in \R Q|\<\b,Q\>\subset \Z\}$ is the dual lattice
of $Q.$ Moreover $e^{2\pi i\b\th}$ maps $S(V)_\a$ to $S(V)_{\alpha+\b}.$
\end{lemma}

\pf See the proof of Proposition 9.5.10 in \cite{PS}. \qed

Consider the representation  $S(V)_\a$ of $\tilde{W}.$ By Theorem
7.6 of \cite{GW} (although Theorem 7.6 of \cite{GW} is stated for
semisimple Lie algebras, but the same argument applies to the case of
Heisenberg algebra), there exists a map $\phi\in \Diff(S^1)\mapsto
\pi_\a(\sigma(\phi))\in U(S(V)_\a)$ which is a unitary cocycle
representation of $\Diff(S^1)$ on $S(V)_\alpha,$ and 
$$\pi_\a(\sigma(\phi))\pi_\a(f(\cdot))\pi_\a(\sigma(\phi)^*)=\pi_\a(f(\phi^{-1}(\cdot))).$$
Let $\B_Q$ be a net on $S(V)_0$ such that
$$\B_Q(I)=\{\pi_0(f)|f\in \tilde{W}, \supp f\subset I\}''$$ where $\pi_0$
denotes the representation of $\tilde{W}$ on  $S(V)_0$.
\begin{proposition}\label{B_Q} We have

(1) $\B_Q$ is a conformal net on $S(V)_0;$ 

(2) $\B_Q$ is strongly additive.
\end{proposition}

\pf (1) is obvious and (2) follows from the same argument of Proposition 1.3.2 of \cite{TL}. \qed

\begin{definition}\label{dA_Q} Let $\A_Q$ be a net of von Neumann algebra on
$H_0$ such that $\A_Q(I)=\{\pi_0(e^{2\pi i f})|e^{2\pi if}\in \L T, \supp f\subset I\}$ where $\pi$ denotes the representation of $\L T$ on $H_0$.
\end{definition}
We note that by definition $\A_Q$ is independent of the choices of
2-cocycle $\epsilon.$\par
By the statement before Proposition \ref{B_Q}, we have a unitary
cocycle representation of $\Diff(S^1)$ on $H_0=\sum_{\a\in Q}S(V)_\a$ such that
$$ \pi(\sigma(\phi))\pi(f)\pi (\sigma(\phi)^*)=\pi(f^{\phi})$$ 
for $f\in \tilde{W}$ where $f^{\phi}(\th):=f({\phi}^{-1}(\th)).$  
We claim that 
$$\pi(\sigma(\phi))\pi(e^{2\pi if})\pi (\sigma(\phi)^*)=c(f,\phi)\pi(e^{2\p if^{\phi}})$$ 
for the some phase factor $c(f,\phi)\in\Co.$  
First we have 
\begin{lemma}\label{locality2} $\A_Q$ is a local net.
\end{lemma}

\pf This follows from Proposition \ref{locality}. \qed

\begin{proposition}\label{locality3} If $\phi\in \Diff(I),$ and $\supp e^{2\pi if}\cap I=\emptyset,$ then 
 $$\pi(\sigma(\phi))\pi(e^{2\pi if})\pi (\sigma(\phi)^*)=\pi(e^{2\p if^{\phi}}).$$ 
\end{proposition}

\pf By Proposition \ref{B_Q}, $\B_Q$ is a conformal net, and it follows that $\pi_0(\sigma(\phi))\in \B_Q(I).$ Note that $\pi$ is a representation of $\B_Q$ on 
$\bigoplus_{\a\in Q}S(V)_\a,$ and restrict to an irreducible representation
of $\B_Q$ on each $S(V)_\a.$ It follows that
$$Ad\pi(\pi_0( \sigma(\phi)))=Ad\pi(\sigma(\phi)).$$
So we have 
$$Ad\pi(\sigma(\phi))(\pi(e^{2\pi if}))=Ad\pi(\pi_0( \sigma(\phi)))(\pi(e^{2\pi if}))=\pi(e^{2\pi if})$$
where we have used Lemma \ref{locality2} in the last equality since
$$\pi_0( \sigma(\phi))\in \B_Q(I)\subset \A_Q(I).$$
\qed

\begin{proposition}\label{Diffcov} We have  for $f\in \L T$ 
$$\pi(\sigma(\phi))\pi(e^{2\pi if})\pi(\sigma(\phi)^*)=c(f,\phi)\pi(e^{2\pi 
i f^\phi})$$
where $c(f,\phi)\in\Co.$
\end{proposition}

\pf Since $\Diff(S^1)$ is a simple group, it is generated by $\Diff(I).$ It is 
sufficient to prove the proposition for $\phi\in\Diff(I).$ Note that 
$$\pi(e^{2\pi if})=\pi(e^{2\pi i(\a\th+f_0)})\pi(e^{2\pi if_1(\th)})
=\pi(e^{2\pi i\a\th})\pi(e^{2\pi i(f_0+f_1)})e^{-\frac{i}{2}\<\a|f_0\>}$$
where $\alpha$ is the winding number of $f$.
Since 
$$Ad\pi(\sigma(\phi))(\pi(e^{2\pi i(f_0+f_1)}))=\pi(e^{2\pi i(f_0+f_1^\phi)})),$$
it is enough to show the proposition for $e^{2\pi i f}=e^{2\pi i\a\th}.$ By the
same argument it reduces to show the proposition for
$e^{2\pi i f}$ so that $f$ has winding number $\alpha$ and 
$\supp e^{2\pi i f}\cap I=\emptyset,$  which follows from Proposition \ref{locality3}. \qed

\begin{proposition}\label{A_Q} (1) $\A_Q$ is a conformal net.

(2) $\A_Q$ is strongly additive and split.
\end{proposition}

\pf (1) follows from Propositions \ref{Diffcov} and \ref{locality3}.
As for (2), let $I_1,I_2$ be two subintervals of $I$ obtained from $I$ 
by removing an interior point of $I.$ By Proposition \ref{B_Q},
$\B(I)=\B(I_1)\vee \B(I_2).$ Since $\A(I)$ is generated by $\B(I)$ and 
$\pi(e^{2\pi if(\theta)})$ with $\supp e^{2\pi if(\theta)}\subset I_1,$ it follows that 
$\A(I_1)\vee \B(I)=\A(I),$ and so $\A(I_1)\vee \A(I_2)=\A(I).$

The character $\Tr q^{L_0}$ of $H_0$ ($L_0$ represents the generator of rotation group $S^1$) is well known to be
$$\Tr q^{L_0}=\frac{\theta_Q(q)}{\eta(q)^l}$$
where 
$$\theta_Q(q)=\sum_{\alpha\in Q}q^{\frac{\<\alpha|\alpha\>}{2}}$$
is the  theta function of the lattice of $Q$ and
$$\eta(q)=q^{\frac{1}{24}}\prod_{n\geq 1}(1-q^n)$$
is the eta function. Here $l$ is the rank of $Q$ and $q=e^{2\pi i\tau}$ 
for $\tau$ in the upper half plane. Hence $q^{L_0}$ is of trace class, and it follows from that $\A_Q$ is split (cf. \cite{BDL}). \qed

\begin{definition}\label{autoLT} 
 Let $\l(\th):[0,1]\to \R Q$ be a smooth
map with $\l(0)=0,$ $\l(1)=\l\in Q^*,$ and $\l^{(n)}(0)=\l^{(n)}(1)=0$ for
all positive $n.$ Define an automorphism
$Ad_{\l(\th)}$ of $\L T$ by the following formula: 
$$Ad_{\l(\th)}(e^{2\pi if},c)=(e^{2\pi if},ce^{2\pi i\int\<\l'|f\>d\th}).$$
\end{definition}

Note that if $\l\in Q^*,$  $e^{2\pi i\l}$ lies in the center of $\L T.$ For any interval $I\subset S^1$ we choose an element $P_{\l,I}\in LT$ such that $P_{\l,I}(\th)=e^{2\pi i\l(\th)}$
if $\th\in I $ where $\l(\th)$ is as in definition \ref{autoLT}.
\begin{definition}\label{auto A_Q} Let $\l\in Q^*.$ Define an automorphism
$Ad_{\l}$ of $\A_Q$ by 
$$Ad_{\l}(y)=\pi(P_{\l,I}) y \pi(P_{\l,I})^*$$
for $y\in \A_Q(I).$
\end{definition}

\begin{lemma}\label{auto} $Ad_\l$ in Definition \ref{auto A_Q} is independent of
the choice of $P_{\l,I}$ and $Ad_{\l}\pi(\L_IT)=\pi(Ad_\l\L_IT)$ for
any $I.$
\end{lemma}

\pf If $P_{\l,I}'$ is another choice, then $P_{\l,I}P_{\l.I}'\in \L_{I'}T$ and
by locality
$$\pi(P_{\l,I}')y\pi(P_{\l,I}')^*=\pi(P_{\l,I})y\pi(P_{\l,I})^*$$ for
$y\in \A_Q(I).$ The equality in the proposition can be checked directly 
by definitions. \qed

\begin{proposition}\label{all irrep of A_Q} 

(1) $Ad_{\l}, \l\in Q^*$ gives an irreducible DHR representation of $\A_Q,$ 
and each such representation corresponds to an irreducible representation
of $\L T,$ labeled by $\l\in Q^*/Q$ (We identify $\l$ with its image in
the quotient map $Q^*\rightarrow  Q^*/Q$ ) 
as in (3) of Lemma \ref{irrep of LT}. 

(2) Let $\psi$ be an irreducible representation of $\A_Q$ with positive 
energy on $H.$ Then $\psi$ is isomorphic to $Ad_{\l}$ for some
$\l\in Q^*/Q.$
\end{proposition}

\pf (1) Note that $e^{2\pi i\l}$ is in the center of $\L T,$ and by the same
proof of Proposition 5.8 of \cite{KLX}, $Ad_{\l}$ is an irreducible DHR
representation of $\A_Q.$ Such a representation of $\A_Q$ corresponds to representation $\pi(Ad_{\l}\L T)$ of $\L T$ by Lemma \ref{auto}.

(2) We will use an idea for the rank one case given in \S4 of
\cite{Xu3}.

Let $f=f_0+f_1(\th): S^1\to \R Q$ be a map 
with winding number $0$. Let $\{I_1,...,I_k\}$ be a finite open
covering of $S^1.$ Assume that $\{\phi_i\}$ is a partition of unity such that
$\supp \phi_i\subset I_i.$ Then $f=\sum_{i=1}^kf\phi_i.$ Let 
$c(f,\phi)\in S^1$ 
be the phase factor in the center of $\L T$ 
so that $(e^{2\pi if},1)=\prod_{j=1}^k(e^{2\pi if\phi_j},1)c(f,\phi)$ as an element of $\L T.$ By using Isotony  we claim 
that the following map 
$$f\to \psi(f)=\prod_{j=1}^k\psi_{I_j}(\pi_{I_j}(e^{2\pi if\phi_j}))
c(f,\phi)$$
is independent of the choice of $\{I_1,...,I_k\}$ and $\{\phi_i\}.$ Moreover,
$$\psi(f)\psi(g)=c(f,g)\psi(f+g).$$
In fact, if $\{J_1,...,J_n\}$ is another open covering of $S^1$ and 
$\{\bar\phi_j\}$ is another partition of unity with $\supp\bar\phi_j
\subset  J_j$ for $j=1,...,n,$ so that $J_{j_1}\cup J_{j_2}
\ne S^1$ for any $1\leq j_1,j_2\leq n.$
We have by Isotony 
$$\pi_{I_s}(e^{2\pi if\phi_s})=\prod_{j=1}^n\pi_{I_s\cap J_j} (e^{2\pi if\phi_s\bar\phi_j})x_i$$
where the phase factor $x_i\in S^1$ is determined by
$$(e^{2\pi if\phi_s},1)=\prod_{j=1}^n( e^{2\pi if\phi_s\bar\phi_j},1)x_i,$$
and $\pi_{I_s\cap J_j} (e^{2\pi if\phi_s\bar\phi_j})\in \A_Q(I_s\cap J_j)$
by definition.
It follows that 
\begin{eqnarray*}
&&\ \ \  \prod_{s=1}^k\psi_{I_s}(\pi_{I_s}(e^{2\pi if\phi_s}))c(f,\phi)\\
&&=\prod_{s=1}^k\prod_{j=1}^n\psi_{I_s\cap J_j}(\pi_{I_s\cap J_j}
(e^{2\pi if\phi_s\bar\phi_j}))x
\end{eqnarray*}
with $x\in S^1.$ Similarly,
\begin{eqnarray*}
&&\ \ \  \prod_{j=1}^n\psi_{J_j}(\pi_{J_j}(e^{2\pi if\bar \phi_j}))c(f,\bar\phi)\\
&&=\prod_{s=1}^k\prod_{j=1}^n\psi_{I_s\cap J_j}(\pi_{I_s\cap J_j}
(e^{2\pi if\phi_s\bar\phi_j}))y
\end{eqnarray*}
for some $y\in S^1.$ Note that $J_{j_1}\cup J_{j_2}
\ne S^1$ for any $1\leq j_1,j_2\leq n.$ It follows that in permuting 
the factors
$$\psi_{I_s\cap J_j}(\pi_{I_s\cap J_j}
(e^{2\pi if\phi_s\bar\phi_j}))$$
above, the phase factors are determined by the group law of $\L T.$ 
Since 
$$\prod_{s=1}^k\prod_{j=1}^n(e^{2\pi if\phi_s\bar\phi_j}, 1)x=
(e^{2\pi if},1)=\prod_{s=1}^k\prod_{j=1}^n(e^{2\pi if\phi_s\bar\phi_j}, 1)y,$$
it follows that
$x=y$ and
$$\prod_{s=1}^k\psi_{I_s}(\pi_{I_s}(e^{2\pi if\phi_s})c(f,\phi)=\prod_{j=1}^n\psi_{J_j}(\pi_{J_j}(e^{2\pi if\bar \phi_j})c(f,\bar\phi).$$
By the independence of $\psi(f)$ on $\{\phi_i\}$ above, it is straightforward
to check that $\psi(f)\psi(g)=c(f,g)\psi(f+g)$ where
$$(e^{2\pi if},1)(e^{2\pi ig},1)=(e^{2\pi i(f+g)},c(f,g))$$
and $\psi(R_{\th})\psi(f) \psi(R_{\th})^*=\psi(f^{\th})$ where
$R_{\th}$ is the rotation by angle $2\pi \th.$

Now for $f=\a\th+f_0+f_1(\th)$  where $\a$ is the winding number of $f$, 
let $g_I(\th)=\a\th+g_0+g_1(\th)$ be such that
$\supp e^{2\pi ig_I}\subset I$ and define
$$\psi(f)=\psi_I(\pi_I(e^{2\pi ig_I}))\psi(f-g_I)c(f,g_I)$$
where $c(f,g_I)$ in the center of $\L T$ is determined by
$$(e^{2\pi if},1)=(e^{2\pi ig_I},1)(e^{2\pi i(f-g_I)},1)c(f,g_I)$$
in $\L T$. 
Note that $f-g_I\in C^{\infty}(S^1,\R Q)$ and $\psi(f-g_I)$ is well-defined
as in the previous paragraph. One checks that $\psi(f)$ is independent
of $I$ and the choice of $g_I,$ and
$$\psi(f)\psi(g)=c(f,g)\psi(f+g),$$
$$\psi(R_{\th})\psi(f) \psi(R_{\th})^*=\psi(f^{\th}).$$

To show that the map $e^{2\pi if}\in \L T\to \psi(f)$ is well-defined
we need to check that if $f=\a,$ then $\psi(\a)$ is the identity operator.
For any $f$ we get
$$\psi(R_{1})\psi(f) \psi(R_{1})^{-1}=\psi(f-\alpha)=\psi(f)e^{2\pi i\frac{1}{2}\<\a,\a\>}\psi(\a).$$
But $\psi(R_1)$ is a scalar operator as $\psi$ is an irreducible representation
of $\A_Q.$ It follows that $\psi(f)=\psi(f)\psi(\a),$ and
$\psi(\a)$ is the identity operator since $\psi(f)$ is unitary. So we
get a well -defined irreducible representation
$$e^{2\pi if}\in \L T\to \psi(f)$$
of $\L T$ with positive energy. By Lemma \ref{irrep of LT} and (1), (2) is
proved. \qed

\begin{remark} There is a similar result in the theory of vertex operator
algebra. Let $V_Q$ be the vertex operator algebra associated to
the lattice $Q$ (cf. \cite{B},\cite{FLM1}). It has been proved in \cite{D}
that $V_{Q+\l}$ for $\l\in Q^*/Q$ gives a complete list of irreducible
$V_Q$-modules up to isomorphism. 
\end{remark}

\begin{lemma}\label{factorize} Let $\A$ be a net with split property. If $\A$ has only finitely many irreducible representations with positive energy
up to isomorphism, then every irreducible representation of $\A\otimes\A$ with positive energy $\pi$ is of the form $\pi_1\otimes\pi_2$ where $\p_i$ 
are irreducible representation of $\A$ with positive energy.
\end{lemma}

\pf As in Lemma 2.7 of \cite{KLM}, it is enough to show that
$\pi(\A\otimes 1)''$ is a type $I$ factor.  Since $\pi_1=\pi\res_{\A\otimes
1}$ is also a representation of $\A$ with positive energy, and $\A$ is split,
by Proposition 5.6 of \cite{KLM} and Lemma 5.14 of \cite{BCL},
$\pi_1=\int_X^{\oplus}\pi_{\l}d\mu(\l)$ where $\pi_{\l}$ are  irreducible representations 
of $\A$ with positive energy for almost all $\l.$ Since $\A$ has only finitely
many irreducible representations with positive energy, it follows that
$\pi(\A\otimes 1)''$ is a type I factor. \qed

\begin{theorem}\label{rational} Let $\A$ be a conformal net, and assume that
$\A$ is strongly additive and split. If $\A$ has  only finitely many irreducible representations with positive energy
up to isomorphism, and if each such representation has finite index, then $\A$ is completely rational, and $\mu_{\A}=\sum_{\lambda}d(\lambda)^2,$ where the sum is over all
irreducible representations of $\A$ with   positive energy.
\end{theorem}

\pf 
The theorem and its proof are essentially contained in \S4 of \cite{LX}
except for the positive energy condition. We will give a proof with 
necessary modifications compared to \cite{LX}. \par
It is sufficient to show that $\mu_{\A}$ is finite. Consider $(\A\otimes \A)^{\Z_2}\subset \A\otimes \A$ where  $(\A\otimes \A)^{\Z_2}$ is the fixed point
subset of  $\A\otimes \A$ under flip. By Corollary 4.6 of \cite{LX}, we have
a representation of  $(\A\otimes \A)^{\Z_2}$ which is a direct sum of two
irreducible representations with positive energy. Let $\tau$  be one 
of them. Note that $\mu_{\A}$ is finite if and only if $\tau$ has finite
index.

By Corollary 3.3 of \cite{LX}, $\a_{\tau^2}$ is a DHR representation of
$\A\otimes\A.$ Since $\pi(\Diff(I))\subset (\A\otimes \A)^{\Z_2}(I),$
by Lemma 4 of \cite{AFK}, $\alpha_{\tau^2}$ is M\"obius invariant with
representation $\pi({\bold G})\in \a_{\tau^2}\res_{(\A\otimes \A)^{\Z_2}}.$ Since
$\a_{\tau^2}|_{(\A\otimes \A)^{\Z_2}}$ is a direct sum of $\tau^2$ and
$\sigma\tau^2$ (cf. \S3.2 of \cite{LX} for the definition
of $\sigma$),  and both have positive energy by Theorem 5.13 of
\cite{BCL}, it follows that $\alpha_{\tau^2}$ is a M\"obius covariant
representation of $\A\otimes A$ with positive energy. Hence
$\a_{\tau^2}=\int_{X}^{\oplus}\l d\mu(\l)$ where almost all $\l$ are
irreducible representations of $\A\otimes \A$ with positive energy. By
Lemma \ref{factorize}, it follows that $\a_{\tau^2}$ is a direct sum of 
irreducible representations of $\A\otimes \A$ with finite index.

Let $[\alpha_{\tau}]=\sum_{i=1}^2[X_i]$ with $X_i$ irreducible as sectors.
Since $[\a_{\tau^2}]= [\alpha_{\tau}][\alpha_{\tau}],$ it follows that
$X_iX_j$ and $X_jX_i$ contain sectors with finite index. By Lemma 
3.6 of \cite{LX}, $X_i$ must have finite index, and so does $\alpha_{\tau}.$
This proves that $\mu_{\A}$ is finite. \qed. 

\begin{corollary}\label{A_Q rational} $\A_Q$ is completely rational and 
$$\mu_{\A_Q}=|Q^*/Q|.$$
\end{corollary}

\pf The proof follows by Propositions \ref{A_Q},\ref{all irrep of A_Q} and
Theorem \ref{rational}.
\qed

\begin{remark} It will be interesting if one can give a direct proof of 
Corollary \ref{A_Q rational} without using Theorem \ref{rational}.
\end{remark}

\section{Orbifolds}

Let $\A_Q$ be the conformal net on $H=\bigoplus_{\a\in Q}S(V)_{\a}$ 
as in Section 3. We will consider a finite automorphism group $\Gamma$ of
$\A_Q$ which arises from isometries of the lattice $Q$ as follows: for each
$\sigma\in\Gamma,$ there is an isometry of $Q$ defined by the same letter $\sigma,$ and
moreover the following map 
$$Ad_\sigma(e^{2\pi i(\Delta_f\th+f_0)} e^{2\pi if_1(\th)}, x)=(\eta(\sigma)^{-1}e^{2\pi i(\sigma(\Delta_f)\th+\sigma(f_0))}e^{2\pi i\sigma(f_1(\th))}, x)$$
gives an automorphism of $\L T$ with finite order. Here $\eta(\a)=\pm 1.$
For such $\sigma\in \Gamma,$ let $\pi(g)$ be the unique unitary operator on $H_0$ such that $\pi(g)\cdot \Omega=\Omega$ ($\Omega$ is the vacuum vector) and 
$$\pi(\sigma)\pi(e^{2\pi if})\pi(\sigma)^*=\pi(Ad_\sigma e^{2\pi if}).$$
One check easily that such unitary operator exists and is unique.   

\begin{lemma}\label{proper action} The map $g\mapsto Ad\pi(g)(y),$ $y\in \A_Q(I)$
defines a proper action of $\Gamma$ on $\A_Q.$ The fixed point subset $\A_Q^{\Gamma}$ is a conformal net.
\end{lemma}

\pf It is enough to show that $\pi(g)\pi(\phi)\pi(g)^*=\pi(\phi)$
for all $\phi\in \Diff(S^1).$ Since $\Diff(I)$ for any $I$ generates $\Diff(S^1)$
and $\Diff(S^1)$ is a perfect group, it suffices to check that
 $\pi(g)\pi(\phi)\pi(g)^*\pi(\phi)^*\in \Co 1$ for all $\phi\in\Diff(S^1).$ Note that  $\pi(g)\pi(\phi)\pi(g)^*\pi(\phi)^*\in \A_Q(I)$ and by definition,
 $\pi(g)\pi(\phi)\pi(g)^*\pi(\phi)^*\in \B_Q(I)'.$  
But  $\B_Q(I)'\cap \A_Q(I)=\Co 1,$ we immediately have 
 $\pi(g)\pi(\phi)\pi(g)^*\pi(\phi)^*\in \Co 1,$ as desired.
\qed
\par
Fix $\sigma\in \Gamma$ of order $N$. Let $P_0:=\frac{1}{N}\sum_{1\leq i\leq N}
\sigma^i$ be the projection on $\R Q.$  For any $\delta\in\R Q$, let
$\delta_*$ be the unique element in the orthogonal complement of 
$P_0(\R Q)$ such that $\delta= \delta_0 + (1-\sigma) \delta_*, 
\delta_0 \in P_0(\R Q).$ We will use $(Q^*/Q)^\sigma$ to denote those 
elements of $Q^*/Q$ which is fixed by $\sigma.$ The set  $(Q^*/Q)^\sigma$
can be represented as follows: let $Q_\sigma:= \{\delta \in Q^*|
(1-\sigma)\delta \in Q \},$ then $(Q^*/Q)^\sigma= Q_\sigma/Q.$
The following lemma follows directly from definitions:
\begin{lemma}\label{id}
$$
(Q^*/Q)/ ( Q_\sigma^*/Q) \simeq Q^*/Q_\sigma^*\simeq  Q_\sigma/Q
= (Q^*/Q)^\sigma
$$
\end{lemma}

\begin{definition}\label{twisted loops} Let $N$ be the order of $\sigma.$ 
Set 
$$C_tT=\{e^{2\pi if}\in C(S^1,T)|f(\th+\frac{1}{N})-\sigma(f(\th))\in Q,  \theta\in[0,1]\}$$
$$C_NT=\{e^{2\pi if}\in C(S^1,T) |f(\th+\frac{1}{N})-f(\th)\in Q, \theta\in[0,1]\}$$
$$C_NT_0=\{e^{2\pi if}\in C_NT|f(0)- \sigma(f(0))\in Q \}$$
For each $e^{2\pi if}\in C_NT_0$ we define a map $\phi_2: C_NT_0\to C_tT$ by
$\phi_2(e^{2\pi if})=e^{2\pi ih}$ where $h$ is the unique continuous function 
which satisfies $h(\th)=f(\th)$ for
$0\leq \th\leq \frac{1}{N}$ and $h(\th+\frac{1}{N})\equiv g(f(\th))$ modulo
$Q$  if $\frac{1}{N}\leq \th\leq 1.$ 
\end{definition}

\begin{example}\label{twist} (1) Let $e^{2\pi iN\a}\in C_tT_0.$ Then $\phi_2(e^{2\pi iN\a})=e^{2\pi ih(\th)}$ such that $h(\theta)=N\a\th$ for $0\leq \th\leq \frac{1}{N}$
and $h(\th+\frac{1}{N})\equiv  \sigma(f(\th))$ modulo $Q$ if
$\frac{1}{N}\leq \th\leq 1.$ There is a unique choice of continuous function $h(\th)$ with the specified property:
$$h(\theta)=\a+\sigma(\a)+\cdots \sigma^{i-1}(\a)+(N\th-i)\sigma^i(\a)$$
for $\frac{i}{N}\leq \th\leq \frac{i+1}{N}$ and $i=0,...,N-1.$
From this example one can see that in general $\phi_2$ maps 
smooth functions to piece-wise smooth functions.\par 

(2) Let  $e^{2\pi ih}\in C_NT_0$ be a constant loop. Then 
$h\equiv g(h)$  modulo $Q$ and $\phi_2(e^{2\pi ih})=e^{2\pi ih}.$
\end{example}
As in section \ref{solitonsection}, let $\xi\in S^1$ and identify
$\R\simeq S^1 \setminus \{\xi\}\simeq (0,1).$

Let
 $$L_{\R}T=\{e^{2\pi if}\in LT|{\supp e^{2\pi if}}\subset (0,1), 
f^{(n)}(0)=f^{(n)}(1)=0, \forall n\geq 1 \}. 
$$
Fix $P=2NQ\subset Q$ which inherits inner product $\<\cdot\>$ from $Q.$ Note that $\<\a,\b\>\in 4N\Z$ for $\a,\b\in P.$  We can choose an equivalence 
class of 2-cocylces $\epsilon$ on $Q$ so that 
$\epsilon (\alpha,\beta)=1, \forall \alpha,\beta\in P$
since $\langle \alpha,\beta\rangle \in 4\Z,  \forall \alpha,\beta\in P.$
Similarly 
for the central extension $\L T$ associated to $P, \langle, \rangle$
(resp.  $\L T$  associated to $P,  \frac{1}{N}\langle, \rangle$), we 
will choose the 2-cocyles as in Definition \ref{MultiplicationofLT} 
to be trivial
since $\langle \alpha,\beta\rangle \in 4N \Z,  \forall \alpha,\beta\in P.$

\begin{definition}\label{Lt}

Denote by $\L_{\R}T$ be the subgroup of
$\L T$ associated to $P, \langle, \rangle$ (with the 
trivial 2-cocycle as above ) whose projection onto $LT$ is
$L_{\R}T.$ Denote by  $\L_N T_0,$ $\L_NT, \L_t T$ the subgroups
of $\L T$ associated with $P, \frac{1}{N} \langle, \rangle$ 
(with the 
trivial 2-cocycle as above )
whose 
projections onto $LT$ are the smooth loops in $C_NT_0,C_NT, C_t T$
respectively.
\end{definition}  

\begin{proposition}\label{N-fold} The homomorphism
$$\begin{array}{ccccc}
L_{\R}T & \stackrel{\phi_1}{\to} & L_N T_0 &\stackrel{\phi_2}{\to} & L_tT \\
e^{2\pi if(\th)}& \stackrel{\phi_1}{\mapsto} & e^{2\pi f(N\th)}& 
\stackrel{\phi_2}{\mapsto} & \phi_(e^{2\pi if(\th)}):=\phi_2( e^{2\pi if(N\th)})
\end{array}$$
can be lifted uniquely to a homomorphism between the central extensions
$$\L_{\R}T  \stackrel{\phi}{\to} \L_tT .$$
\end{proposition}

\pf
Note that $\phi_2\phi_1$ maps $L_\R T$ to smooth loops in
$LT.$
The proof is a direct computation by definitions. \qed

Let $\A_P$ be the conformal subnet of $\A_Q$ associated to the lattice $P$
and represented on its vacuum Hilbert space $H_P,$  and $\A_{P,\frac{1}{N}}$
the conformal net associated to $(P,\frac{1}{N}\<\cdot\>)$ on
$H$ where the extra subscript $\frac{1}{N}$ indicates that the
inner product on $P$ is $\frac{1}{N}\<\cdot\>.$  Using $\phi_1$ we get a covariant representation $\pi$ of $\A_P$ 
on $H$ such  that 
$$\pi (\pi_P(e^{2\pi if(\th)}))=\pi(\phi_1(e^{2\pi if}))=\pi(e^{2\pi if(N\th)})$$
for any localized $e^{2\pi if}.$ Denote by $Pr$  the projection from $H$ to
$\overline{\pi(\L_tT)\Omega}.$
\begin{definition}\label{basic soliton} Let $e^{2\pi if}\in \L_{\R}T,$ and
$$\pi(e^{2\pi if(N\th)})=\prod_{i=0}^{N-1}Ad_{R^i}(\pi(e^{2\pi i\hat f(\th)}))$$
with $\hat f(\th)=f(N\th)$ if $0\leq \th\leq \frac{1}{N},$ and 
$\hat f(\th)=\Delta_f$ if $\frac{1}{N}\leq \th\leq 1,$ and $R$ is the rotation by $\frac{2\pi }{N}.$ Define
$$\pi_\sigma^{(\xi)}(e^{2\pi if(\th)})=\prod_{i=0}^{N-1}Ad_{\sigma^iR^i}(\pi(e^{2\pi i \hat f(\th)}))Pr.$$

\end{definition}

\begin{proposition}\label{Normality1} $\pi_\sigma^{(\xi)}$ extends to a normal
representation of $\A_P(\R)$ and restricts to a DHR representation of $\A_P^{\Gamma}.$
\end{proposition}

\pf Fix $I\subset \R\simeq(0,1)\simeq S^1\setminus \{\xi\},$ 
and assume that $\supp e^{2\pi if(\th)}\subset I.$ Then 
$$\pi^{(\xi)}(e^{2\pi if(N\th)})=\prod_{i=0}^{N-1}Ad_{R^i}(\pi(e^{2\pi i\hat f(\th)}))
\subset \A_{P,\frac{1}{N}}(I_1)\vee \cdots \vee  \A_{P,\frac{1}{N}}(I_N)$$
where $I_i^N=I,$ $i=1,...,N,$ as intervals on $S^1,$ and we choose the ordering so that on $(0,1),$ $I_i$ is to the left of $I_{i+1}.$ Since  
$\A_{P,\frac{1}{N}}$ is split by Prop. \ref{A_Q}, 
there is a normal isomorphism 
$$\chi_I:  \A_{P,\frac{1}{N}}(I_1)\otimes\cdots \otimes  \A_{P,\frac{1}{N}}(I_N) \to  \A_{P,\frac{1}{N}}(I_1)\vee\cdots\vee  \A_{P,\frac{1}{N}}(I_N)$$
so that $\chi_I(x_1\otimes\cdots \otimes x_N)=x_1\cdots x_N, 
x_i\in  \A_{P,\frac{1}{N}}(I_i), i=1,...,N. $ Note that
the map $U_{\sigma}:   x_1\otimes\cdots \otimes x_N \to 
x_1\otimes \sigma(x_2)\otimes \cdots \otimes \sigma^{N-1}(x_N)$ is normal. So we have 
$$\pi_\sigma^{(\xi)}(e^{2\pi if(\th)})=\chi_IU_{\sigma}\chi_I^{-1}(\pi(e^{2\pi if(N\th)}))$$
if $\supp e^{2\pi if}\subset I.$ Thus $\pi_{\sigma}^{(\xi)}(x)=\chi_IU_{\sigma}\chi_I^{-1}
(\pi(\phi_1(x)))$ is a normal representation of $\A_P(I).$ 

Note that 
$$\prod_{i=0}^{N-1}Ad_{\sigma^iR^i}(Ad_{R^j}x)=\prod_{i=0}^{N-1}Ad_{\sigma^iR^i}(Ad_{\sigma^{-j}}x)$$
for $x\in \A_{P,\frac{1}{N}}(I_1)\vee\cdots\vee  \A_{P,\frac{1}{N}}(I_N)$ since $\sigma$ and $R$ commute. 

Also note that $Ad_{\sigma}\pi(\phi_1(x))=\pi(\phi_1(Ad_{\sigma}x))$ for 
$x\in \A_P(I).$ So if $x\in \A_P^{\<\sigma\>}(I),$ then 
$$\prod_{i=0}^{N-1}Ad_{\sigma^iR^i}(Ad_{R^j}(\pi(\phi_1(x))))
=\prod_{i=0}^{N-1}Ad_{\sigma^iR^i}(\pi(\phi_1(x)))).$$
So when we change $\xi\in S^1$ to $\xi^1\in S^1$ we have $\pi_\sigma^{(\xi)}=\pi_\sigma^{(\xi^1)}$ when restricting  to $ \A_P^{\<\sigma\>}.$ It follows that 
 $\pi_\sigma^{(\xi)}$ restricts to a DHR representation of  $ \A_P^{\<\sigma\>}.$
\qed

\begin{remark} Note that the definition of soliton above is similar to the soliton given in \cite{LX},  but without using $\Diff(S^1)$ and hence are 
different. We will write    $\pi_\sigma^{(\xi)}$ simply as $\pi_\sigma$ when restricting to  $ \A_P^{\<\sigma\>}.$
\end{remark}

\begin{proposition}\label{Irred} $\pi_\sigma^{(\xi)}$ is an irreducible soliton
 of $\A_P.$ 
\end{proposition}

\pf First we prove that the representation $\pi(\L_t T)Pr$ on $PrH$ is irreducible. Notice that this is a representation with positive energy, and the 
identity component 
$(\L_t T)^o$
is the product of a  Heisenberg group and  $T_0:= \exp^{2\pi i P_0(\R P)}.$ 
It follows from (1)-(2) of
Lemma \ref{irrep of LT} that the representation
$\pi (\L_t T)Pr$ of $\L_t T$ is irreducible. To prove the proposition it is enough to show that $\pi_\sigma^{(\xi)}(\A_P(\R))''=\pi(\L_tT)''Pr,$
and it is sufficient to show that $\pi(e^{2\pi ig})Pr \in\pi_\sigma^{(\xi)}(\A_P(\R))''$ for any $e^{2\pi ig}\in (\L_tT)^o.$ 

Let $x(\th)$ be a complex valued function on $[0,1]$ with $x(\th+\frac{1}{N})=x(\th),$ and $x^{(n)}(0)=x^{(n)}(\frac{1}{N})=0$ for all $n.$ It follows 
from Lemma 1.2.2 of \cite{TL} that for any $\inz>0$ one can choose 
$x_{\inz}$ such that $||x_{\inz}-1||_{\frac{1}{2}}<\inz$ 
(cf. \S1.2 of  \cite{TL} for the definition of norm $||.||_{\frac{1}{2}}$),
and by Proposition 1.3.2 of \cite{TL} we have that
$\pi(e^{2\pi ix_{\inz}g})\to \pi(e^{2\pi ig})$ strongly. Note that 
 $\pi(e^{2\pi ix_{\inz}g})=\pi (\phi_2\phi_1(e^{2\pi if_{\inz}})),$ 
where $e^{2\pi if_{\inz}}\in \L_{\R}T$ with $f_{\inz}(\th)=x_{\inz}(\frac{\th}{N})g(\frac{\th}{N})$ for $0\leq \th\leq 1.$ It follows that 
$\pi(e^{2\pi ig})Pr\in \pi_\sigma^{(\xi)}(\A_P(\R))''.$ \qed
\begin{definition}\label{Auto L_tT} Let $\l(\th):[0,1]\to \Co P$ be a smooth
map with $\l(0)=0,$ $\l(1)=\l\in P^*,$ and $\l^{(n)}(0)=\l^{(n)}(1)=0$ for
all positive $n.$ Define an automorphism $Ad_{\l(\th)}$ on $\L_tT$ such that
$$Ad_{\l(\th)}(e^{2\pi i(NP_0(\a)\th+\a_*)},x)=(e^{2\pi i(NP_0(\a)\th+\a_*)},xe^{2\pi i\int\<\l'|P_0(\a)\th+\a_*\>d\th}),$$
$$Ad_{\l(\th)}(e^{2\pi ih(\th)},y)=(e^{2\pi ih(\th)},ye^{2\pi iN\int_{0}^{\frac{1}{N}}\<\l'(N\th)|h(\th)\>d\th})$$
where  $h(\th+\frac{1}{N})=\sigma h(\th)$ for $0\leq \th\leq 1,$ $\int hd\th=0.$
\end{definition}

\begin{lemma}\label{Autocov} For any $e^{2\pi if}\in \L_{\R}T$ we have 
$$\pi(\phi_2\phi_1(Ad_{\l(\th)}e^{2\pi if}))=\pi(Ad_{\l(\th)}\phi_2\phi_1
( e^{2\pi if})).$$
\end{lemma}

\pf It is straightforward by definition. \qed
\begin{proposition}\label{Irred2}  $\pi_\sigma^{(\xi)}(Ad_{\l})$ is an irreducible soliton
 of $\A_P(\R)$ which corresponds to an irreducible representation of $\L_tT$ 
where the central group $(P^*/P)^{\sigma}$ of $\L_tT$ acts as 
$\mu \in  (Q^*/Q)^{\sigma}\mapsto e^{2\pi i\<\mu|\l\>}.$
\end{proposition}

\pf  $\pi_\sigma^{(\xi)}(Ad_{\l})$  is irreducible since  $\pi_\sigma^{(\xi)}$ is by Proposition \ref{Irred}. The second statement follows by
lemma \ref{Autocov}.
\qed

To make contact with the results in Section 4 of \cite{BK}, and to motivate
the definitions in the next section, let us define (compare to 
Section 4.3 of \cite{BK}).

\begin{definition}\label{G_P} Let $G_P=S^1\times \exp(2\pi iP_0(\R P))\times P$ be the set consisting of elements $ce^hU_{\alpha}$ $(c\in S^1, h\in 2\p i P_0(\R P),
\a\in P).$ Define a multiplication in $G_P$ by the formulas
$$e^{h}e^{h'}=e^{h+h'}$$
$$e^hU_{\a}e^{-h}=e^{-\<h|\a\>}U_{\a}$$
$$U_{\a}U_{\b}=U_{\a+\b}.$$
\end{definition}
Note that $G_P$ is a group. 
Recall the Heisenberg group $\tilde V_Q$ associated with a lattice $Q$
and inner product $\langle,\rangle$ on $Q$ before Lemma \ref{irrep of LT}. 
We will consider a  Heisenberg group $\tilde{V}_{P,\frac{1}{N}}$
associated with the lattice $P$ and inner product $\frac{1}{N} \langle,\rangle$ on $P.$
\begin{lemma}\label{iso} We have the isomorphism
$$\L_tT\cong G_P\times \tilde{V}_{P,\frac{1}{N}}.$$
\end{lemma} 

\pf 
First we note that 
$\L_tT$ is generated by $e^{2\pi i(NP_0(\a)\th+\a_*)}$ and $e^{2\pi ih(\th)}$ 
with $h(\th+\frac{1}{N})=\sigma h(\th)$ for $0\leq \th\leq 1,$ $\int hd\th=0.$
For any $\mu \in  (P^*/P)^{\sigma},$ $e^{2\pi i\mu}$ is
in the center of $\L_tT.$

Let $\psi:\L_tT\to   G_P\times \tilde{V}_{P,\frac{1}{N}}$ be a map such that
$$\phi(e^{2\pi ih},x)=(e^{2\pi ih},x)\in \tilde{V}_{P,\frac{1}{N}}$$
$$\phi(e^{2\pi i(N\pi_0(\a)\th+\a_*)}, y)=yU_{\a}\in G_P$$
$$ \phi(e^{2\pi ih_0})=e^{2\pi ih_0}\in G_P$$
for $h(\th)$ with $h(\th+\frac{1}{N})=\sigma h(\th), 0\leq \th\leq 1,$ 
$\int hd\th=0,$ and $h_0\in P_0(\R P).$ One checks directly by definition that $\phi$ is
an isomorphism of groups. \qed

We note that a class of irreducible representations of  $G_P\times 
\tilde{V}_{P,\frac{1}{N}}$ is given by Theorem 4.2 of \cite{BK}, and these irreducible
representations are determined uniquely by the action of central subgroup 
$(Q^*/Q)^{\sigma}$ of $G_P$ as follows: Given $\mu \in Q^*/Q,$ there
is an irreducible representation $\pi_{\mu}$ of $G_P\times \tilde{V}_{P,\frac{1}{N}}$ on a Hilbert space $H_{\mu}$ such that $\l \in (Q^*/Q)^{\sigma}$
acts as $e^{2\pi i\<\l|\mu\>}.$ By using Proposition \ref{Irred2} and Lemma
\ref{iso}, it is easy to check that $\pi_{\mu}$ corresponds to $\pi_{\sigma}(Ad_{\mu(\th)}).$ Denote by $\pi_{\sigma,\mu}= \pi_{\sigma}(Ad_{\mu(\th)})$
and note that using this notation $\pi_{\sigma}=\pi_{\sigma,0}.$ 

\begin{proposition}\label{sigma_1} Let $\sigma_1\in\Gamma.$ Then 

(1) $\pi_{\sigma,\mu}(Ad\sigma_1)\cong\pi_{\sigma_1^{-1}\sigma\sigma_1,\sigma_1^{-1}\mu}$ as solitons of $\A_P(\R).$

(2) $\pi_{\sigma,\l_1}\cong\pi_{\sigma,\l_2}$ if and only if $\l_1-\l_2\in P_{\sigma}^*.$
\end{proposition}

\pf (1) First note that $Ad_{\mu(\th)}Ad_{\sigma_1}=Ad_{\sigma_1}Ad_{\sigma_1^{-1}\mu(\th)}$ on $\L_{\R}T,$ and so it is sufficient to show that
$\pi_{\sigma}(Ad_{\sigma_1})\cong\pi_{\sigma_1^{-1}\sigma\sigma_1}.$ Fix $I\subset \R,$ let $I_i$ be intervals on $S^1$ so that $I_i^N=I,$ $i=1,...,N.$ Then 
$$\prod_{i=1}^NAd_{\sigma^iR^i}(Ad_{\sigma_1}x)=Ad_{\sigma_1}(\prod_{i=1}^NAd_{(\sigma_1^{-1}\sigma\sigma_1)^iR^i}(x))$$ 
for all $x\in A_{P,\frac{1}{N}}(I_1\vee\cdots\vee I_N).$ It follows that $\pi_{\sigma}(Ad_{\sigma_1})\cong\pi_{\sigma_1^{-1}\sigma\sigma_1}.$

(2) By Proposition \ref{Irred2}, $\pi_{\sigma,\l_1}\cong\pi_{\sigma,\l_2}$
if and only if  $\<\mu|\l_1-\l_2\>\in \Z$ for all $\mu\in (P^*/P)^{\sigma}
=P_{\sigma}/P.$ It follows that $\pi_{\sigma,\l_1}\cong\pi_{\sigma,\l_2}$
if and only if $\l_1-\l_2\in P_{\sigma}^*.$ 
\qed

Consider $\A_P^{\<\sigma\>}\subset \A_P,$ the fixed point subnet of $\A_P$
under the action of $\sigma.$ Let $V\in \A_P(I)$ 
be a unitary such that $\sigma(V)=e^{\frac{2\pi i }{N}}V$ (cf. 
\S8 of \cite{LX}).

$Ad_V$ induces a DHR representation of  $\A_P^{\<\sigma\>}.$
By Proposition \ref{Normality1} and Prop. 8.2 
of \cite{LX} $\pi_\sigma$ decomposes into a direct sum of
$N$ irreducible representations of $\A_P^{\<\sigma\>}.$ Let
$\pi_{\sigma0}$ be one of such irreducible representations of $\A_P^{\<\sigma\>} $ obtained from $\pi_\sigma$ by projecting onto $\sigma$ invariant
subspace.   
Set $\pi_{\nu}=\pi_{\sigma0}(Ad V).$ By Proposition \ref{Normality1}, $\pi_{\nu}$
restricts to a DHR representation of  $\A_P^{\<\sigma\>},$ and
it is a covariant representation of  $\A_P^{\<\sigma\>}.$ 
In fact, $\pi(\phi(\cdot))$ is a covariant representation of  
$\A_P^{\<\sigma\>}$
 as $\A_P^{\<\sigma\>}$ is conformal  and   $\pi(\phi(\cdot))$ restricts
to a DHR  representation of  $\A_P^{\<\sigma\>}$ (cf. \cite{AFK}).

 Let
$$U_I=\{g\in {\bold G}|gI\cup I\subset S^1-\{\xi\}\}$$
be a neighborhood of identity in ${\bold G}.$ For $g\in U_I$ we have 
$$\pi_{\nu}(gxg^*)=\pi_{\nu}(g)\pi_{\nu}(x)\pi_{\nu}(g^*)=
\pi_{\sigma 0}(g)\pi_{\sigma 0}(g^*VgV^*)\pi_{\nu}(x)\pi_{\sigma 0}(Vg^*V^*g)\pi_{\sigma 0}(g)^*$$
for all $x\in  \A_P^{\<\sigma\>}.$ So we have 
$ \pi_{\nu}(g)^*\pi_{\sigma 0}(g)\pi_{\sigma 0}(g^*VgV^*)^*\in \Co 1$ for all $g\in U_I.$
It follows that $ \pi_{\nu}(g)=\pi_{\sigma 0}(g)\pi_{\sigma 0}(g^*VgV^*)$ for 
all $g\in U_I$ as the only one dimensional representation of $\bold {G}$ is
the trivial representation.

Consider $\th\to \pi_{\sigma 0}(R_{\th})$ 
and set 
$$F(\th)= Ad_{\pi(r_{-\frac{\th}{N}})}(\pi(\phi(V)))\pi(\phi(V))^*Pr_{\sigma0}$$
where $\pi(r_{-\frac{\th}{N}})$ denotes the unitary operator on $H$ 
implementing rotation by $-\frac{\th}{N}$ on $\A_{P,\frac{1}{N}},$  
and $Pr_{\sigma0}=Pr \sum_{1\leq \sigma \leq N} \sigma^i$ is the projection 
onto the irreducible representation $\pi_{\sigma 0}$ of $\A_P^{\<\sigma\>}.$ 
Note that
if $\th\in U_I$ then $F(\th)=\pi_\sigma(R_{\th})\pi_{\sigma 0}(R^*_{\th}V R_{\th}V^*).$ 
Also $\pi(\phi(AdR_{\th},y))=\pi(Ad_{r_{\frac{\th}{N}}}(\phi(y)))$ for all
$y\in \A_P^{\<\sigma\>}$ by definitions, i.e.,
$$\pi(R_{\th})\pi(\phi(y))\pi(R_{\th})^*=\pi(r_{\frac{\th}{N}}) \pi(\phi(y))\pi(r_{\frac{\th}{N}})^*$$ for $y\in \A_P^{\<\sigma\>}.$ Using this one checks easily that 
$$F(\th_1+\th_2)=F(\th_1)F(\th_2).$$
It follows that $F(\th)=\pi_{\nu}(R_{\th}),$ for all $\th,$ since both sides
are one-parameter group of unitaries which agree on a neighborhood of $0.$

On the other hand,
\begin{eqnarray*}
&&F(2\pi)=\pi_{\sigma 0}(R_{2\pi}Ad_{\pi(r_{-\frac{2\pi}{N}})}(\pi(\phi(V)))
\pi(\phi(V^*))\\
&&\ \ \ \ \  \ =\pi_{\sigma 0}(R_{2\pi}\pi(\phi(\sigma(V)V^*))\\
&&\ \ \ \ \ \ =e^{\frac{2\pi i}{N}}
\pi_{\sigma 0}(R_{2\pi}).
\end{eqnarray*} 
where we have used $\sigma(V) V^*= e^{\frac{2\pi i}{N}}.$
It follows that the univalence of $\pi_{\sigma 0}(Ad V)$
is the univalence of $\pi_{\sigma 0}$ multiplied by $e^{\frac{2\pi i}{N}}.$ By monodromy equation (cf. \cite{FRS}), we have proved the following 
\begin{proposition}\label{Grading} $G(\pi_{\sigma 0}, AdV)
=e^{\frac{2\pi i}{N}},$ where $G(\cdot,\cdot)$ is defined as in Lemma 8.3 of
\cite{LX}.
\end{proposition}

\begin{remark} The same argument as in the proof of Proposition \ref{Grading}
shows that $K(1)=1$ in the paragraph after (47) of \cite{LX}.
\end{remark}
\subsection{General case}\label{gs}

When $\frac{1}{N}\<\cdot\>$ is not an even integral on $Q,$ we do not have an analogue of Proposition \ref{N-fold} and net $\A_{Q,\frac{1}{N}}.$ However we have a subnet $\A_{2NQ}\subset \A_Q$ where we can apply the construction 
of Section 2.1. Also there is no $\L_tT$ in general case, but we have $\tilde{V}_{Q,\frac{1}{N}}$ and an analogue $G_Q$ as defined as follows (cf. Section
4.3 of \cite{BK}).

\begin{definition}\label{G_Q} Let $G_Q=S^1\times \exp(2\pi iP_0(\R Q))\times Q$ be the set consisting of elements of the form  $ce^hU_{\alpha},$ $(c\in S^1, h\in 2\pi i P_0(\R Q), \a\in P)$ with multiplication
 $$e^{h}e^{h'}=e^{h+h'}$$
$$e^hU_{\a}e^{-h}=e^{-\<h|\a\>}U_{\a}$$
$$U_{\a}U_{\b}=\inz(\a,\b)e^{\pi i\<\a|\b_*\>-\frac{1}{2}\<\a|\b\>+\frac{1}{2}\<\a|P_0(\b)\>}U_{\a+\b}.$$
\end{definition}
One checks easily that $G_Q$ is a group.
\begin{remark} Our group $G_Q$ is slightly different from $G$ of Section
4.3 of \cite{BK}, the commutator among $U_{\a}, U_{\b}$ is the complex conjugate of 4.44 of \cite{BK}. The reason for defining the multiplication rule for $U_{\a}U_{\b}$ comes from the following computations: Let $h(\th)$ be as in 
(1) of Examples \ref{twist}, then $h(\th)=N\pi_0(\a)\th+\a_*+h_{\a}(\th).$ 
Regarding $e^{2\pi ih_{\a}(\th)}$ as an element in $\tilde{V}_{Q,\frac{1}{N}}$ 
we have 
$$e^{2\pi ih_{\a}(\th)}e^{2\pi ih_{\b}(\th)}=e^{\pi i\<\a|\b_*\>-\frac{1}{2}\<\a|\b\>+\frac{1}{2}\<\a|P_0(\b)\>}e^{2\pi i(h_{\a}(\th)+h_{\b}(\th))}.$$ Hence if we map
$e^{i\a\th}$ to $U_{\a}e^{2\pi ih_{\a}(\th)},$ to preserve the commutator relations we need 
$U_{\a}U_{\b}=\inz(\a,\b)e^{\pi i\<\a|\b_*\>-\frac{1}{2}\<\a|\b\>+\frac{1}{2}\<\a|P_0(\b)\>}U_{\a+\b}.$
\end{remark}

We will treat $G_{P}$ (cf, Definition \ref{G_P}) as a subgroup of $G_Q$ under the natural map $G_P\to G_Q.$ 

According to Lemma \ref{iso}, the analogue of $\L_tT$ is now replaced by the group $G_Q\times \tilde{V}_{Q,\frac{1}{N}}.$ 

\begin{definition}\label{Auto G_Q} For $\l(\th)$ as in Definition \ref{Auto L_tT}
with $\l\in Q^*$ define an automorphism on $G_Q\times \tilde{V}_{Q,\frac{1}{N}}$ by
$$Ad_{\l(\th)}U(\a)=U(\a)e^{2\pi i\int\<\l'|P_0(\a)\th+\a_*\>}$$
$$Ad_{\l(\th)}e^{ih}=e^{ih}e^{2\pi iN\int_0^{\frac{1}{N}}\<\l'(N\th)|h(\th)\>d\th}$$
for $h$ with $h(\th+\frac{1}{N})=\sigma(h(\th)).$
\end{definition}

\begin{definition}\label{phi} 
Let $f(\th)=\Delta_f\th+f_0+f_1.$
Define $\hat f(\th)=f(N\th),$ $0\leq \th\leq \frac{1}{N}$
and $\hat f(\th)$ is a continuous function  with $\hat f(\th+\frac{1}{N})=\sigma(\hat f(\th))$
modulo $Q.$ It follows that 
$$\hat f(\th)=NP_0(\a)\th+P_0(f_0)+(\Delta_f)_*+\hat f_1(\th)$$
where  $\hat f_1(\th+\frac{1}{N})=\sigma(\hat f_1(\th))$ and $\int_{S^1} 
\hat f_1d\th=0.$
Then the map $\phi$ is defined as 
$$\phi(e^{if(\th)})=U(\Delta_f)e^{2\pi i(P_0(f_0)+\hat f_1(\th))}
e^{-\pi i \<P_0(\Delta f)|f_0\>}.$$
\end{definition}

\begin{lemma}\label{homo} $\phi:\L_\R T\to  G_Q\times \tilde{V}_{Q,\frac{1}{N}}$
is a group homomorphism and $\phi(Ad_{\l(\th)}x)=Ad_{\l(\th)}\phi(x)$ for
all $x\in \L_\R T.$ 
\end{lemma}

\pf The proof is a direct (but tedious) check by using definitions. \qed

Fix $\sigma\in \Gamma$ with $\sigma^N=1.$ Let $\pi_{\l}$ be an irreducible
representation of $G_Q\times \tilde{V}_{Q,\frac{1}{N}}$ as given by Theorem
4.2 of \cite{BK}. This is an irreducible representation of $ G_Q\times \tilde{V}_{Q,\frac{1}{N}}$ on a Hilbert space $H_{\l}$ where the central subgroup 
$(Q^*/Q)^{\sigma}$ of $G_Q$ acts as the character $e^{2\pi i\<\l|\mu\>}$ for
$\mu\in (Q^*/Q)^{\sigma}.$ This representation, when restricting to 
$G_{P}  \times \tilde{V}_{Q,\frac{1}{N}},$ decomposes into direct sum of finitely many irreducible representations of  $G_{P} \times \tilde{V}_{Q,\frac{1}{N}}:$ 
$$H_{\l}=\bigoplus_{\omega}H_{\l,\omega}\otimes K_{\omega}$$
where $K_{\omega}$ is an irreducible representation of 
$G_P\times \tilde{V}_{Q,\frac{1}{N}}$ as given by 
Theorem 4.2 of \cite{BK} for lattice $P=2NQ,$ and the central subgroup 
$(P^*/P)^{\sigma}$ of $G_P$ acts by the character $e^{2\pi i\<\mu_1|\omega\>}$
for $\mu_1\in (P^*/P)^{\sigma}$ and each $H_{\l,\omega}$ is of finite dimensional. We note that by Proposition \ref{Irred2}, $K_{\omega}$ corresponds to 
representation $\pi_{\sigma,\omega}$ of $\A_P(\R).$ 
Fix an interval $I\subset S^1-\{\xi\},$ and a set of representatives $\a_1,\cdots, \a_k,$ $k=({\rm rank}\,Q)^{2N}$ for the finite abelian group $Q/P.$ 
By abuse of notations, in this section 
we will use $ \L_\R T$ to be the central extension
of $L_\R T$ as in Definition \ref{Lt}, but associated with lattice $Q.$
Choose 
$e^{2\pi if_{\a_i}}\in \L_\R T$ so that $\Delta_{f_{\a_i}}=\a_i,$
$i=1,...,k,$ and $\supp e^{2\pi if_{\a_i}}\subset I.$ Note that for each $I\subset J\subset \R^1-\{\xi\},$ every element $x\in \A_Q(J)$ can be written uniquely as 
$x=\sum_{s=1}^kx_i\pi(e^{2\pi if_s}),$ where $x_i\in \A_P(J).$

\begin{definition}\label{generalsoliton} With notations as above, we define 
$$\pi_{\sigma,\l}^{(\xi)}(x)=\sum_{1\leq s\leq k,\omega}
Id_{\l,\omega}\otimes \pi_{\sigma,\omega}^{(\xi)}(x_i)\pi_{\l}(e^{2\pi if_s})$$ 
for $x\in \A_Q(J).$
\end{definition}

\begin{proposition}\label{generalnormal} $\pi_{\sigma,\l}^{(\xi)}$ 
as defined in 
Definition \ref{generalsoliton} is an irreducible soliton representation
of $\A_Q.$  $\pi_{\sigma,\l}^{(\xi)}
(Ad\mu({\th}))\cong \pi_{\sigma,\l+\mu}^{(\xi)},$
$\pi_{\sigma,\l}^{(\xi)}(Ad\sigma_1)=\pi^{(\xi)}_{\sigma^{-1}_1 \sigma \sigma_1,\sigma^{-1}_1\l}$ for $\sigma_1\in \Gamma,$ 
and $\pi_{\sigma,\l_1}^{(\xi)} \simeq \pi_{\sigma,\l_2}^{(\xi)}$ iff
$\l_1-\l_2\in Q_\sigma^*.$ Moreover $\pi_{\sigma,\l}^{(\xi)}$ 
restricts to a DHR representation of $\A_Q^{\langle\sigma\rangle}.$
\end{proposition}

\pf 
Let $J$ be the interval as defined before Definition \ref{generalsoliton}.  
By Propositions \ref{Normality1} and \ref{Irred2}, the map $x\in \A_{Q}(J)\to
\pi_{\sigma,\l}^{(\xi)}(x)$ is normal. To check that $\pi_{\sigma,\l}^{(\xi)}$ is a
homomorphism, it is enough to check
$\pi_{\sigma,\l}^{(\xi)}(x_1x_2)=\pi_{\sigma,\l}^{(\xi)}(x_1)\pi_{\sigma,\l}^{(\xi)}(x_2)$ for
$x_1,x_2$ in a set of generators for $\A_Q(J).$ We can choose this set
of generators to be elements of $\L_J T.$ It follows from Lemma
\ref{homo} that $\pi_{\sigma,\l}$ is a homomorphism. The rest of the
proposition follows from Propositions \ref{Normality1}, \ref{sigma_1}
and definitions. \qed

We will use $\pi_{\sigma,\l}$ to denote the  DHR representation of 
$\A_Q^{\langle\sigma\rangle}$ from Prop. \ref{generalnormal} and let $\pi_{\sigma0}$ 
be an irreducible subrepresentation of $\pi_{\sigma,0}|_{\A_Q^{\sigma}}.$

\begin{proposition}\label{Grading general} $G(\pi_{\sigma0},AdV)
=e^{\frac{2\pi i}{N}}$  where $G(\cdot,\cdot)$ is defined as in Lemma 8.3 of
\cite{LX}.
\end{proposition}

\pf Since $\pi_{\sigma0}(R_{\th})\in  \pi_{\sigma 0}(\A_Q^{\sigma})'',$ it follows by Proposition \ref{Grading} and definition of $\pi_{\sigma}$ that the univaluence of $\pi_{\sigma 0}( \A_Q^{\sigma})$ is the univalence of $\pi_{\sigma 0}$ multiplied by $e^{\frac{2\pi i}{N}},$ and this implies the proposition by monodromy equation (cf. \cite{FRS}). \qed

\subsection{Irreducible representations of $\A_Q^{\Gamma}.$}\label{generalsolitonsection}
By Prop.\ref{generalnormal},  
$\pi_{\sigma,\l_1}^{(\xi)} \simeq \pi_{\sigma,\l_2}^{(\xi)}$ iff
$\l_1-\l_2\in Q_\sigma^*.$ By using lemma \ref{id}, we will identify
$\l$ with its image in $(Q^*/Q)^\sigma$ under the composition 
of quotient map 
$Q^*\rightarrow Q^*/ Q_\sigma^*$ and 
$Q^*/ Q_\sigma^*\simeq  Q_\sigma /Q=(Q^*/Q)^\sigma$ 
in this section. 

\begin{theorem}\label{identity} Let $\sigma_i\in\Gamma,$ $\l_i\in (Q^*/Q)^{\sigma},$ and $\pi_{\sigma_i,\l_i}^{(\xi)}$ be solitons of $\A_Q^{\Gamma}$ as in proposition
\ref{generalnormal}. Then $\pi_{\sigma_1,\l_1}^{(\xi)}
\cong\pi_{\sigma_2,\l_2}^{(\xi)}$
as solitons of $\A_Q$ if and only if $\sigma_1=\sigma_2$ and $\l_1=\l_2.$
\end{theorem}

\pf The proof is similar to that of Theorem 7.1 of \cite{KLX}. It is sufficient to show that $\pi_{\sigma_1,\l_1}^{(\xi)}\cong\pi_{\sigma_2,\l_2}^{(\xi)}$ implies $\sigma_1=\sigma_2.$ Note that $\pi_{\sigma_i,\l_i}^{(\xi)}$ restrict to DHR representations of 
$\A_Q^{\<\sigma_i\>}.$ It follows that $\pi_{\sigma_1,\l_1}^{(\xi)}$ restricts to a DHR
representation of $\A_Q^{\<\sigma_1\>}\vee \A_Q^{\<\sigma_2\>},$ and 
by Proposition
\ref{Grading general} and the proof (2) of Proposition 7.2 of \cite{KLX} we 
must have  $\A_Q^{\<\sigma_1\>}\vee \A_Q^{\<\sigma_2\>}=\A_Q^{\<\sigma_1\>}.
$ By Galois 
correspondence (cf. \cite{ILP}) we have $\<\sigma_2\>\subset \<\sigma_1\>.$ Exchanging 
$\sigma_1$ and $\sigma_2$ we conclude that $\<\sigma_2\>=\<\sigma_1\>.$ So
$\sigma_2=\sigma_1^m$ for some integer $m$ with $(m,order(\sigma_1))=1.$ By 
Prop. \ref{Grading general}, the same
proof as in (1) of Proposition 7.2 of \cite{KLX} shows that $m=1.$ \qed
\par
We will use the following simple lemma, and we refer the reader to 
\S2 of \cite{KLX} for definitions of endomorphisms.
\begin{lemma}\label{abelian act} Let $\rho\in \End(M)$ be an endomorphism
and assume that $[\rho\tau_i]=[\rho]$ for $i=1,...,k$ where $\tau_i\in \Aut(M)$ and $[\tau_i]\ne [\tau_j]$ if $i\ne j.$ Then $d^2_{\rho}\geq k.$ 
\end{lemma}

\pf If $d_{\rho}=\infty$ there is nothing to prove. Assume that $d_{\rho}$ is finite. By by Frobenius duality $[\bar\rho\rho]$ contains $\bigoplus_{i=1}^k[\tau_i]$ and so     $d^2_{\rho}\geq k.$ \qed

\begin{proposition}\label{index of soliton 1} Let $\pi_{\sigma,\l}^{(\xi)}$ be as in Proposition \ref{generalnormal}. Then $d(\pi_{\sigma,\l}^{(\xi)})^2\geq \frac{|Q*/Q|}{
|(Q^*/Q)^{\sigma}|}.$
\end{proposition}

\pf By Proposition \ref{generalnormal}, 
$$\pi_{\sigma,\l}^{(\xi)}(Ad_{\mu({\th})})\cong \pi_{\sigma,\l+\mu}^{(\xi)}$$ 
for $\mu\in Q^*/Q.$  It follows that $\pi_{\sigma,\l}^{(\xi)}(Ad_{\mu({\th})})\cong \pi_{\sigma,\l}^{(\xi)}$ if and only if $\mu\in Q_{\sigma}^*/Q.$ 
The inequality now follows 
from Lemma \ref{abelian act} and 
Lemma \ref{id}.
\qed

Let $\pi_{\sigma_i,\l_i}^{(\xi)}$ be solitons of $\A_Q$ $(i=1,2).$ By
Proposition 4.4 of \cite{KLX}, $\pi_{\sigma_1,\l_1}^{(\xi)}\res
\A^{\Gamma}\cong \pi_{\sigma_2,\l_2}^{(\xi)}\res \A^{\Gamma}$ if and only if
there exists $\sigma\in \Gamma$ such that
$\pi_{\sigma_1,\l_1}^{(\xi)}(Ad\sigma)=\pi_{\sigma_2,\l_2}^{(\xi)}.$  
By Proposition
\ref{generalnormal} and Theorem \ref{identity} we have
$\sigma^{-1}\sigma_1\sigma=\sigma_2$ and $\sigma^{-1}\l_1=\l_2.$ This
motivates the following consideration. Let
$\{[\sigma_1],...,[\sigma_m]\}$ be a list of conjugacy classes in
$\Gamma.$ Each $\pi_{\sigma_i,\l_i}$ restricts to a DHR representation
of $\A_Q^{\Gamma}.$ Let $\Gamma_{\sigma_i}:= \{ \sigma\in \Gamma |  \sigma\sigma_i=\sigma_i \sigma
\}$ and  let $[\l_{i,1}],...,[\l_{i,m_i}]$ be the orbits of $(Q^*/Q)^{\sigma_i}$
under the action of $\Gamma_{\sigma_i}.$ Let $\Gamma_{\sigma_i,\l_{i,s}}:=\{\sigma\in
\Gamma_{\sigma_i}|\sigma\l_{i,s}=\l_{i,s}\}$ for $s=1,...,m_i.$ Then for any
$\sigma\in \Gamma_{\sigma_i,\l_{i,s}},$ $\pi_{\sigma_i,\l_{i,s}}$ and
$\pi_{\sigma_i,\l_{i,s}}(Ad\sigma)$ are isomorphic as representations
of $A^{\Gamma}_{Q}$ and by Theorem 4.8 of \cite{KLX}, there is a
projective representation of $\Gamma_{\sigma_i,\l_{i,s}}$ on
$H_{\sigma_i,\l_{i,s}}$ with cocycle $c_{\pi_{\sigma_i,\l_{i,s}}}$
such that $H_{\sigma_i,\l_{i,s}}=\bigoplus_{\delta}M_{\delta}\otimes
H_{(\delta,\sigma_i,\l_{i,s})}$ where each $M_{\delta}$ is an irreducible
representation of $\Gamma_{\sigma_i,\l_{i,s}}$ with cocycle
$c_{\pi_{\sigma_i,\l_{i,s}}}$ (and all irreducible representations of
$\Gamma_{\sigma_i,\l_{i,s}}$ with the same  cocycle
$c_{\pi_{\sigma_i,\l_{i,s}}}$ appear in the decomposition of
$H_{\sigma_i,\l_{i,s}}$), and $H_{(\delta,\sigma_i,\l_{i,s})}$ supports an 
irreducible representation of $\A_Q^{\Gamma}$ with 
$H_{(\delta,\sigma_i,\l_{i,s})}\not\cong H_{(\delta',\sigma_i,\l_{i,s})}$
if $\delta\ne \delta'.$ 
By Theorem 4.5 of \cite{KLX} and Prop. \ref{index of soliton 1}, for 
fixed $\sigma_i, \l_{i,s},$ $s=1,...,m_i,$
$$\sum_{\delta}d(\delta,\sigma_i,\l_{i,s})^2=\frac{|\Gamma|^2}{|\Gamma_{\sigma_i,\l_{i,s}}|}d(\pi_{\sigma_i,\l_{i,s}})^2\geq \frac{|\Gamma|^2}{|\Gamma_{\sigma_i,\l_{i,s}}|}\frac{|Q^*/Q|}{|(Q^*/Q)^{\sigma_i}|}.$$
By definition we have
$$|(Q^*/Q)^{\sigma_i}|=\sum_{1\leq s\leq m_i}\frac{|\Gamma_{\sigma_i}|}{|\Gamma_{\sigma_i,\l_i,s}|},$$
and so 
$$\sum_{\delta, 1\leq s\leq m_i}d(\delta,\sigma_i,\l_{i,s})^2\geq \frac{|\Gamma|^2}{|\Gamma_{\sigma_i}|}|Q^*/Q|.$$

Now sum over $1\leq i\leq m$ and note that $|\Gamma|=\sum_{1\leq i\leq m} \frac{|\Gamma|}{|\Gamma_{\sigma_i}|}.$ We have
$$\sum_{\delta, 1\leq i\leq k, 1\leq s\leq m_i}d(\delta,\sigma_i,\l_{i,s})^2\geq |\Gamma|^2|Q^*/Q|=\mu_{\A_Q^{\Gamma}}$$
where we have used Cor. \ref{A_Q rational} and Th. \ref{orb}.
It follows that all $``\geq''$ above are $``=''$ by Theorem \ref{orb},
and  we have 
proved the following theorem:
\begin{theorem}\label{allirrep} $H_{\delta,\sigma_i,\l_{i,s}}$ for 
$i=1,...,m,$ $\l_i\in (Q^*/Q)^{\sigma_i},$ $1\leq s\leq m_i$ as above give 
a complete list of irreducible representations of $\A_Q^{\Gamma},$
and these representations generate a unitary modular tensor 
category. 
\end{theorem}

From the proof and Prop. \ref{index of soliton 1} we have also
proved: 
\begin{corollary}\label{index of soliton 2} $$d(\pi_{\sigma,\l}^{(\xi)})=\sqrt{\frac{|Q^*/Q|}{|(Q^*/Q)^{\sigma}|}}.$$
\end{corollary}

\end{document}